\documentclass[a4paper,leqno,10pt]{article}
\usepackage{amssymb}
\usepackage{amsmath}
\usepackage[german]{babel}
\begin{document}
\def\refname{ }
\begin{center}
\begin{Large}
\textbf{On a result of Fel'dman on linear forms in the values of some $E$-functions}
\end{Large}
\vskip0.4cm
\sc{Keijo V\"a\"an\"anen}
\end{center}
\vskip0.6cm

\begin{center}
Abstract
\end{center}

We shall consider a result of Fel'dman, where a sharp Baker-type lower bound is obtained for linear forms in the values of some $E$-functions. Fel'dman's proof is based on an explicit construction of Pad\'{e} approximations of the first kind for these functions. In the present paper we introduce Pad\'{e} approximations of the second kind for the same functions and use these to obtain a slightly improved version of Fel'dman's result.\\



\section{Introduction}

In 1964 Baker \cite{B} studied linear forms $x_1e^{\alpha_1} + \cdots + x_me^{\alpha_m}$, where $(x_1,\ldots,x_m) \in \mathbb{Z}^m\setminus{\{\underline{0}\}}$ and $\alpha_j \ (j = 1,\ldots,m)$ are distinct rational numbers, and proved a lower bound
\begin{equation}\label{1}
\left|x_1e^{\alpha_1} + \cdots + x_me^{\alpha_m}\right| > h^{1 - c_0/\sqrt{\log\log h}}\prod_{j=1}^m h_j^{-1},
\end{equation}
for all $h = \max\{\left|x_j\right|\} \geq c_1 > e$, where $h_j = \max\{1,\left|x_j\right|\}$ and $c_0, c_1$ are positive constants depending on $\alpha_j$. These constants were made completely explicit in Mahler \cite{M1}. Lower bounds like above depending on each individual coefficient $x_j$ are called Baker-type lower bounds. Baker's proof used essentially Siegel's method with a new idea in the construction of the auxiliary function, a Pad\'{e} type approximation of the first kind for the functions $e^{\alpha_jz}$, obtained by using Siegel's lemma. After that the same idea was used to study other $E$- and $G$-functions satisfying linear differential equations of first order with rational coefficients, see for example \cite{M} and \cite{V}. Then, in an important and deep paper \cite{Z}, Zudilin was able to obtain a similar result for the values of a class of $E$-functions satisfying a system of homogeneous linear differential equations with rational coefficients, in this general result the term $\sqrt{\log\log h}$ in the bound is replaced by $(\log\log h)^{1/(m^2-m+2)}$.

Shortly after Baker's work Fel'dman \cite{F} considered linear forms of the values of the $E$-functions
\begin{equation}\label{2}
\varphi_{\lambda_j}(z) = \sum_{\nu=0}^\infty \frac{z^\nu}{[\nu]_j}, \quad j = 1,\ldots,m,
\end{equation}
where
\[
[0]_j = 1, [\nu]_j = (1+\lambda_j)\cdots(\nu+\lambda_j), \nu \geq 1,
\]
and $\lambda_j \neq -1,-2,\ldots$ are rational numbers such that $\lambda_i - \lambda_j \notin \mathbb{Z}$, if $i \neq j$. Instead of using Siegel's lemma he constructed explicitly appropriate Pad\'{e} approximations of the first kind for the functions $\varphi_{\lambda_j}(z)$ and by using these obtained the following result.\\

\noindent
\textbf{Theorem} (Fel'dman). \textit{Let $\alpha \neq 0$ be a rational number. There exists a positive constant $c_0$ depending on $\lambda_1,\ldots,\lambda_m, m$ and $\alpha$ such that, for all $(x_0,x_1,\ldots,x_m) \in \mathbb{Z}^{m+1}\setminus{\{\underline{0}\}}$, 
\begin{equation}\label{3}
\left|x_0 + x_1\varphi_{\lambda_1}(\alpha) + \cdots + x_m\varphi_{\lambda_m}(\alpha)\right| > H^{-1-c_0/\log\log (H+2)},
\end{equation}
where $H = \prod_{j=1}^m h_j, h_j = \max\{1,\left|x_j\right|\} \ (j=1,\ldots,m)$.}\\

This seems to be still the only result of this type for $E$-functions, where $\sqrt{\log\log}$ in the estimate is improved to $\log\log$. Our main purpose in this paper is to give a new proof for the above Feldman's theorem, where we explicitly construct Pad\'{e} approximations of the second kind for the functions $\varphi_{\lambda_j}(z)$, in other words, simultaneous rational approximations to the functions $\varphi_{\lambda_j}(z)$, which are suitable for proving Baker-type bounds. The application of \cite[Corollary 3.5]{Ma} then leads to a slightly more precise form of the above Theorem, where $c_0$ is given explicitly for large $H$.\\

\noindent
\textbf{Theorem 1}. \textit{Assume that $\lambda_1,\ldots, \lambda_m$ satisfy the assumptions of Fel'dman's theorem. Let $K$ denote $\mathbb{Q}$ or an imaginary quadratic field and $\mathbb{Z}_K$ the ring of integers of $K$, and let $\alpha \in K\setminus \{0\}$. Then there exists a positive constant $H_0$ depending on $\lambda_1,\ldots, \lambda_m, m$ and $\alpha$ such that, for all $(\beta_0,\beta_1,\ldots,\beta_m) \in \mathbb{Z}_K^{m+1}\setminus\{\underline{0}\}$ with $H = \prod_{j=1}^m h_j \geq H_0, h_j = \max \{1,\left|\beta_j\right|\} \ (j=1,\ldots,m)$,
\[
\left|\beta_0+\beta_1\varphi_1(\alpha)+\cdots+\beta_m\varphi_m(\alpha)\right| > H^{-1-\frac{6(d_0+d_1m+d_2m^2)}{\log\log H}},
\]
where $d_0, d_1, d_2$ are positive constants depending on $\lambda_1,\ldots, \lambda_m$ and $\alpha$, to be given explicitly at the end of Section 6.}\\

Pad\'{e} approximations of the second kind were first used in the connection of Baker-type bounds in Sorokin \cite{S} to the consideration of some $G$-functions. Then in \cite{VZ} such a construction was used to study certain $q$-series, for a refinement see also \cite{L}. Moreover, the paper \cite{S} on $\varphi_\lambda(z)$ and \cite{ELM} on the exponential function also apply Pad\'{e} type approximations of the second kind to improve the constants in the above results of Baker and Mahler. In these papers Sorokin used explicit construction but all other applied Siegel's lemma. In fact, as far as we know, the explicit construction of the approximations of the second kind below is the first one for Baker-type bounds of $E$-functions.\\

\section{Explicit construction 1}

Let $n_1,\ldots,n_m$ denote positive integers, $N = n_1+\cdots+n_m$, and 
\[
Q_0(z) = \sum_{k=0}^N a_kz^k.
\]
By denoting $\varphi_j(z) = \varphi_{\lambda_j}(z)$ we have
\[
Q_0(z)\varphi_j(z) = \sum_{\mu=0}^\infty c_{j\mu}z^\mu, \quad c_{j\mu} = \sum_{k=0}^{\min\{\mu,N\}} \frac{a_k}{[\mu-k]_j}, \ j = 1,\ldots,m.
\]
To get the needed Pad\'{e} approximations of the second kind we now choose the coefficients $a_k$ in such a way that $c_{j\mu} = 0$ for all $\mu = N+1,\ldots,N+n_j, j = 1,\ldots,m$. This means that 
\[
a_0 + a_1(\mu+\lambda_j) + a_2(\mu+\lambda_j)(\mu+\lambda_j-1) +\cdots + a_N(\mu+\lambda_j)\cdots(\mu+\lambda_j-(N-1)) = 0
\]
for all $\mu = N+1,\ldots,N+n_j, j = 1,\ldots,m$. This is a system of $N$ linear homogeneous equations in $N+1$ unknowns $a_k$, which has a non-trivial solution. To determine such a solution we denote
\[
\gamma_1 = N+1+\lambda_1, \ldots, \gamma_{n_1} = N+n_1+\lambda_1,
\]
\[
\gamma_{n_1+1} = N+1+\lambda_2, \ldots, \gamma_{n_1+n_2} = N+n_2+\lambda_2, \ldots
\]
\[
\gamma_{n_1+\cdots+n_{m-1}+1} = N+1+\lambda_m, \ldots, \gamma_N = N+n_m+\lambda_m.
\]
Then the above system of equations can be given in the form
\begin{equation}\label{4}
a_0 + a_1\gamma_i + a_2\gamma_i(\gamma_i-1) +\cdots+ a_{N-1}\gamma_i\cdots(\gamma_i-(N-2)) = -a_N\gamma_i\cdots(\gamma_i-(N-1)), \ i = 1,\ldots,N.
\end{equation}
The coefficient determinant $\delta$ of this system is
\[
\delta = \det(1 \ \gamma_i \ \gamma_i(\gamma_i-1) \ldots \gamma_i\cdots(\gamma_i-(N-2)))_{i=1,\ldots,N} = \prod_{1\leq i<j \leq N}(\gamma_j-\gamma_i) \neq 0.
\]
After the choice of $a_N$ we thus obtain a unique solution $a_0,a_1,\ldots,a_{N-1}$.
 
For $\sigma = 1,\ldots,N$, let $\delta_\sigma(z)$ denote the determinant obtained from $\delta$ after replacing $\gamma_\sigma$ by $z$. Then
\[
\delta_\sigma(z) = \delta_{\sigma0} + \delta_{\sigma1}z + \delta_{\sigma2}z(z-1) + \cdots + \delta_{\sigma,N-1}z(z-1)\cdots(z-(N-2)),
\]
where $\delta_{\sigma k}$ is the cofactor of $\delta$ corresponding to the $\sigma,k$-entry ($\sigma = 1,\ldots,N; k = 0,\ldots,N-1$). Since $\delta_\sigma(\gamma_s) = 0$ for all $s \neq \sigma$, we have
\[
\delta_\sigma(z) = c \prod_{s=1,s\neq \sigma}^N(z - \gamma_s)
\]
with some constant $c$, and since $\delta = \delta_\sigma(\gamma_\sigma)$,
\[
c = \delta \prod_{s=1,s\neq \sigma}^N(\gamma_\sigma - \gamma_s)^{-1}.
\]
Thus we get
\begin{equation}\label{5}
\delta_{\sigma0} + \delta_{\sigma1}z + \delta_{\sigma2}z(z-1) + \cdots + \delta_{\sigma,N-1}z(z-1)\cdots(z-(N-2)) = \delta \prod_{s=1,s\neq \sigma}^N \frac{z - \gamma_s}{\gamma_\sigma - \gamma_s}.
\end{equation}

By choosing $z = \kappa$ in (\ref{5}) for each $\kappa = 0,1,\ldots,N-1$, we obtain  
\[
\delta_{\sigma0} + \kappa\delta_{\sigma1} + \kappa(\kappa-1)\delta_{\sigma2} +\cdots+ \kappa!\delta_{\sigma \kappa} = \delta \prod_{s=1,s\neq \sigma}^N \frac{\kappa - \gamma_s}{\gamma_\sigma - \gamma_s}.
\]
So
\begin{equation}\label{6}
A(\frac{\delta_{\sigma0}}{\delta},\frac{\delta_{\sigma1}}{\delta},\ldots,\frac{\delta_{\sigma,N-1}}{\delta})^T = 
\end{equation}
\[
(\frac{1}{0!}\prod_{s=1,s\neq \sigma}^N \frac{- \gamma_s}{\gamma_\sigma - \gamma_s},\frac{1}{1!}\prod_{s=1,s\neq \sigma}^N \frac{1 - \gamma_s}{\gamma_\sigma - \gamma_s},\ldots,\frac{1}{(N-1)!}\prod_{s=1,s\neq \sigma}^N \frac{N -1 - \gamma_s}{\gamma_\sigma - \gamma_s})^T,
\]
where $A$ is the $N\times N$-matrix with rows
\[
(\frac{1}{\kappa!}, \frac{1}{(\kappa-1)!}, \cdots, \frac{1}{1!}, \frac{1}{0!}, 0,\ldots, 0), \ \kappa=0,1,\ldots,N-1.
\]
We now see that $A^{-1}$ is the matrix with rows
\[
((-1)^k\frac{1}{k!},(-1)^{k-1}\frac{1}{(k-1)!}, \cdots, -\frac{1}{1!}, \frac{1}{0!}, 0,\ldots, 0), \ k=0,1,\ldots,N-1,
\]
and therefore the above equality (\ref{6}) implies
\begin{equation}\label{7}
\frac{k!\delta_{\sigma k}}{\delta} = \sum_{\tau = 0}^k(-1)^{k-\tau}{k \choose \tau}\prod_{s=1,s\neq \sigma}^N \frac{\tau - \gamma_s}{\gamma_\sigma - \gamma_s}, \quad k = 0,1,\ldots,N-1.
\end{equation}

By using Cramer's rule we obtain from (\ref{4})
\[
a_k = -a_N\sum_{\sigma =1}^N\frac{\delta_{\sigma k}}{\delta}\prod_{\mu =0}^{N-1}(\gamma_\sigma - \mu), \ k = 0,1,\ldots,N-1.
\]
The choice $a_N = -1/N!$ together with (\ref{7}) then gives, for all $k = 0,1,\ldots,N-1$,
\begin{equation}\label{8}
k!a_k = \sum_{\sigma =1}^N \sum_{\tau = 0}^k(-1)^{k-\tau}{k \choose \tau}\prod_{\mu =0}^{N-1}\frac{\gamma_\sigma - \mu}{1+\mu}\prod_{s=1,s\neq \sigma}^N \frac{\tau - \gamma_s}{\gamma_\sigma - \gamma_i}.
\end{equation}
Thus we have explicitly constructed polynomials
\[
Q_0(z) = \sum_{k=0}^N a_kz^k, \ P_{0j}(z) = \sum_{\mu=0}^N c_{j\mu}z^\mu, \ j = 1,\ldots,m,
\]
such that $\deg Q_0(z) = N, \deg P_{0j}(z) \leq N$, and the remainder terms
\[
R_{0j}(z) := Q_0(z)\varphi_j(z) - P_{0j}(z) = \sum_{\mu=N+n_j+1}^\infty c_{j\mu}z^\mu, \ j = 1,\ldots,m.\\
\]

\section{Explicit construction 2}

The construction above is not enough, since we need $m+1$ linearly independent approximations. To get these we fix $i, 1 \leq i \leq m$, and denote
\[
\gamma_0 = N+1+\lambda_i,
\]
\[
\gamma_1 = N+1+\lambda_1+\delta_{1i}, \ldots, \gamma_{n_1} = N+n_1+\lambda_1+\delta_{1i},
\]
\[
\gamma_{n_1+1} = N+1+\lambda_2+\delta_{2i}, \ldots, \gamma_{n_1+n_2} = N+n_2+\lambda_2+\delta_{2i}, \ldots
\]
\[
\gamma_{n_1+\cdots+n_{m-1}+1} = N+1+\lambda_m+\delta_{mi}, \ldots, \gamma_N = N+n_m+\lambda_m+\delta_{mi},
\]
where $\delta_{ij}$ denotes Kronecker's $\delta$. Instead of (\ref{4}) we now consider the system of equations
\[
a_0 + a_1\gamma_0 + a_2\gamma_0(\gamma_0-1) +\cdots+ a_N\gamma_0(\gamma_0-1)\cdots(\gamma_0-(N-1)) = 1,
\]
\[  
a_0 + a_1\gamma_\sigma + a_2\gamma_\sigma(\gamma_\sigma-1) +\cdots+ a_N\gamma_\sigma(\gamma_\sigma-1)\cdots(\gamma_\sigma-(N-1)) = 0, \ \sigma = 1,\ldots,N,
\]
with a coefficient determinant
\[
\Delta = \prod_{0\leq \ell<j \leq N}(\gamma_j-\gamma_\ell) \neq 0.
\]
By Cramer's rule this system has a solution
\[
a_k = \frac{\Delta_{0k}}{\Delta}, \quad k = 0,1,\ldots,N,
\]
where $\Delta_{0k}$ is the cofactor of $\Delta$ corresponding to the $0,k$-entry. To give $a_k$ explicitly we proceed as in the previous section. Analogously to (\ref{5}) we now have
\[
\Delta_{00} + \Delta_{01}z + \Delta_{02}z(z-1) +\cdots+ \Delta_{0N}z(z-1)\cdots(z-(N-1)) = \Delta \prod_{s=1}^N \frac{z-\gamma_s}{\gamma_0-\gamma_s}.
\]
Repeating the considerations leading to (\ref{7}) we then obtain
\begin{equation}\label{9}
k!a_k = \sum_{\tau = 0}^k(-1)^{k-\tau}{k \choose \tau}\prod_{s=1}^N \frac{\tau - \gamma_s}{\gamma_0 - \gamma_s}, \quad k = 0,1,\ldots,N.
\end{equation}

For each $i = 1,\ldots,N$ we have thus constructed polynomials
\[
Q_i(z) = \sum_{k=0}^N a_{ik}z^k, \ P_{ij}(z) = \sum_{\mu=0}^{N+\delta_{ij}} c_{ij\mu}z^\mu, \ j = 1,\ldots,m,
\]
where $a_{ik} = a_k$ are given in (\ref{9}), $c_{ij\mu} = c_{j\mu}$ (with $a_k = a_{ik}$), $\deg Q_i(z) = N \ (a_N = 0$ implies $a_0 = \cdots = a_{N-1} = 0$), $\deg P_{ii}(z) = N+1, \deg P_{ij}(z) \leq N$ for all $i \neq j$, and the remainder terms
\[
R_{ij}(z) := Q_i(z)\varphi_j(z) - P_{ij}(z) = \sum_{\mu=N+n_j+1+\delta_{ij}}^\infty c_{ij\mu}z^\mu, \ j = 1,\ldots,m.\\
\]
These approximations and the approximation of the previous section satisfy the following lemma.\\

\noindent
\textbf{Lemma 1}. \textit{The determinant}
\[
\Omega(z) = \det (Q_i(z) \ P_{i1}(z) \ \ldots \ P_{im}(z))_{i=0,1,\ldots,m} = cz^{(m+1)N+m},
\]
where
\[
c = \frac{-1}{N!}\prod_{i=1}^m\prod_{\nu=1}^{N+1}(\lambda_i+\nu)^{-1}.\\
\]

\textit{Proof}. The coefficients of the leading terms of $Q_0(z)$ and $P_{ii}(z) \ (i = 1,\ldots,m)$ are $-1/N!$ and $1/((\lambda_i+1)\cdots(\lambda_i+N+1))$, respectively, here we use the first equation above satisfied by $a_{ik}$. Therefore $\Omega(z)$ is a polynomial of exact degree $(m+1)N+m$ and the coefficient of the leading term of $\Omega(z)$ is the product of the above coefficients.

On the other hand
\[
\Omega(z) = (-1)^m \det (Q_i(z) \ R_{i1}(z) \ \ldots \ R_{im}(z))_{i=0,1,\ldots,m}.
\]
Since ord $R_{ij}(z) \geq N+n_j+1$ and $N = n_1+\cdots+n_m$, it follows that ord $\Omega(z) \geq (m+1)N+m$. This proves Lemma 1.\\

\section{Denominators and upper bounds}

We first give a lemma from \cite [pp. 145-147]{M2} considering the quotients
\[
\frac{(\alpha+1)_n}{n!} =: \frac{u_n}{v_n}, \ (u_n,v_n) = 1, v_n \geq 1, n = 0,1,\ldots,
\]
where $\alpha = r/s \neq -1,-2,\ldots$ with integers $r$ and $s \geq 1, (r,s) = 1$, and $(\alpha)_0 = 1, (\alpha)_n = \alpha(\alpha+1)\cdots(\alpha+n-1)$ for $n \geq 1$.\\

\noindent
\textbf{Lemma 2}. \textit{Let
\[
U_n = \prod_{p\nmid s}p^{[\log(\left|r\right|+sn)/\log p]}, \ V_n = s^{2n}.
\]
Then the least common multiples of $u_0,u_1,\ldots,u_n$ and of $v_0,v_1,\ldots,v_n$ are divisors of $U_n$ and $V_n$, respectively.}\\ 

Let us denote
\[
\lambda_j = \frac{r_j}{s_j}, (r_j,s_j) = 1, s_j \geq 1, \quad \lambda_k-\lambda_j = \frac{r_{kj}}{s_{kj}}, (r_{kj},s_{kj}) = 1, s_{kj} \geq 2.
\]
Further, let
\[
R = \max \{\left|r_j\right|\}, S = \max \{s_j\}, \quad \hat{R} = \max \{\left|r_{kj}\right|\}, \hat{S} = \max \{s_{kj}\}.
\]
Clearly $\hat{R} \leq 2RS$ and $\hat{S} \leq S^2$.

We now consider the denominators of $k!a_k = k!a_{ik}$ in (\ref{9}). Here the product
\[
\Pi_{i\tau} := \prod_{s=1}^N \frac{\tau - \gamma_s}{\gamma_0 - \gamma_s} = \prod_{j=1}^m\prod_{\nu=1}^{n_j}\frac{N+\nu+\lambda_j+\delta_{ji}-\tau}{\lambda_j-\lambda_i+\nu+\delta_{ji}} = 
\]
\[
\prod_{j=1,j\neq i}^m(\frac{(\lambda_j+N+1-\tau)_{n_j}}{n_j!}\cdot\frac{n_j!}{(\lambda_j-\lambda_i+1)_{n_j}})\cdot \frac{(\lambda_i+N+1-\tau)_{n_i+1}}{(n_i+1)!}.
\]
By Lemma 2, the denominator of $\Pi_{i\tau}$ is a factor of
\[
(\prod_{j=1}^m s_j^{2(n_j+\delta_{ji})})\cdot\prod_{j=1,j\neq i}^m\prod_{p} p^{[\log(\left|r_{ji}\right|+s_{ji}n_j)/\log p]}.
\]
Thus the denominators of all $\Pi_{i\tau}$ are factors of 
\begin{equation}\label{10}
D_1 := \prod_{j=1}^m (s_j^{2(n_j+1)}\prod_{p} p^{[\log(\hat{R}+\hat{S}n_j)/\log p]}),
\end{equation}
and so, by (\ref{9}), all $k!D_1a_{ik} \in \mathbb{Z} \ (k=0,1,\ldots,N; i =1,\ldots,m)$. By the weak form of the prime number theorem, see for example \cite[p. 296]{Bu}, the number of primes $p \leq x$
\[
\pi(x) \leq 8\log 2\frac{x}{\log x} < \frac{6x}{\log x}
\]
for all $x > 1$, and therefore
\begin{equation}\label{11}
D_1 \leq S^{2(N+m)}e^{6(\hat{R}m+\hat{S}N)} =: E_1.
\end{equation}

By Lemma 2 and the above expression for $\Pi_{i\tau}$ we also have
\[
\left|\Pi_{i\tau}\right| \leq \prod_{j=1,j\neq i}^m(\frac{s_{ji}^{2n_j}}{s_j^{n_j}}\prod_{p} p^{[\log(R+S(N+n_j))/\log p]})\cdot \frac{1}{s_i^{n_i}}\prod_{p} p^{[\log(R+S(N+1+n_i))/\log p]} 
\]
\[
\leq S^{3N}e^{6(Rm+S+S(m+1)N)}.
\]
This implies, by (\ref{9}),
\begin{equation}\label{12}
\left|k!a_{ik}\right| \leq 2^kS^{3N}e^{6(Rm+S+S(m+1)N)} =: 2^kF_1, \quad k=0,1,\ldots,N; i =1,\ldots,m,
\end{equation}
and so
\begin{equation}\label{13}
\left|Q_i(z)\right| \leq \sum_{k=0}^N \left|a_{ik}z^k\right| \leq F_1e^{2\left|z\right|}.
\end{equation}

Next we consider the coefficients of the polynomials $P_{ij}(z)$,
\[
c_{ij\mu} = \sum_{k=0}^\mu \frac{a_{ik}}{[\mu-k]_j} = \sum_{k=0}^\mu (\frac{k!a_{ik}}{k!(\mu-k)!}\cdot\frac{(\mu-k)!}{(\lambda_j+1)\cdots(\lambda_j+\mu-k)}), \ \mu=0,1,\ldots,N,
\]
remember also, that $c_{ii,N+1} = 1/(\lambda_i+1)_{N+1}$. By Lemma 2 and the above considerations
\[
(N+1)!D_2c_{ii,N+1}, \ N!D_2c_{ij\mu} \in \mathbb{Z}, \quad 1 \leq i,j \leq m; \mu = 0,1,\ldots,N,
\]
where
\begin{equation}\label{14}
D_2 := D_1\prod_{p} p^{[\log(R+S(N+1))/\log p]} \leq S^{2(N+m)}e^{6(R+S+\hat{R}m+(\hat{S}+S)N)} =: E_2, 
\end{equation}
to get this upper bound we used (\ref{11}). Thus
\[
(N+1)!D_2Q_i(z), \ (N+1)!D_2P_{ij}(z) \in \mathbb{Z}[z], \quad i, j = 1,\ldots,m.
\]

Finally we need to consider the polynomials $Q_0(z)$ and $P_{0j}(z)$ constructed in Section 2, here the coefficients $a_k$ are given in (\ref{8}). If $\gamma_\sigma = N+\kappa+\lambda_t, 1 \leq \kappa \leq n_t$, then the last product in (\ref{8}) is
\[
\Pi_{\sigma\tau}^* = \prod_{s=1,s\neq \sigma}^N \frac{\tau - \gamma_s}{\gamma_\sigma - \gamma_s} = \prod_{j=1,j\neq t}^m\prod_{\nu=1}^{n_j}(\frac{N+\nu+\lambda_j-\tau}{\lambda_j-\lambda_t+\nu-\kappa})\cdot\frac{(-1)^{\kappa-1}}{N+\kappa+\lambda_t-\tau}\cdot\frac{n_t!}{(\kappa-1)!(n_t-\kappa)!}\cdot 
\]
\[
\prod_{\nu=1}^{n_t}\frac{N+\nu+\lambda_t-\tau}{\nu} = \prod_{j=1,j\neq t}^m(\frac{(\lambda_j+(N+1-\tau)_{n_j}}{n_j!}\cdot\frac{n_j!}{(\lambda_j-\lambda_t+1-\kappa)_{n_j}})\cdot\frac{(-1)^{\kappa-1}s_t}{r_t+(N+\kappa-\tau)s_t}\cdot
\]
\[
\frac{n_t!}{(\kappa-1)!(n_t-\kappa)!}\cdot \frac{(\lambda_t+N+1-\tau)_{n_t}}{n_t!}.
\]
Since, for all $t, \kappa$ and $\tau$, the number $r_t+(N+\kappa-\tau)s_t$ is a factor of
\[
\prod_{p}p^{[\log (R+2NS)/\log p]},
\]
it follows by Lemma 2 that the denominators of all $\Pi_{\sigma\tau}^*$ are factors of
\[
\prod_{j=1}^m (s_j^{2n_j}\prod_{p} p^{[\log(\hat{R}+\hat{S}N)/\log p]})\prod_{p}p^{[\log (R+2NS)/\log p]}. 
\]
Moreover
\[
\prod_{\mu =0}^{N-1}\frac{\gamma_\sigma - \mu}{1+\mu} = \frac{(\lambda_t+1+\kappa)_N}{N!},
\]
and so Lemma 2 and (\ref{8}) imply that all $k!D_1^*a_{k} \in \mathbb{Z} \ (k=0,1,\ldots,N)$, where
\begin{equation}\label{15}
D_1^* = (s_1\cdots s_m)^{2N}\prod_{j=1}^m (s_j^{2n_j}\prod_{p} p^{[\log(\hat{R}+\hat{S}N)/\log p]})\prod_{p}p^{[\log(R+2NS)/\log p]} 
\end{equation}
\[
\leq S^{2(m+1)N}e^{6(m(\hat{R}+\hat{S}N)+R+2SN)} =: E_1^*.
\]
Note here, that $D_1 \mid D_1^*$.

We now use once again Lemma 2 to get
\[
\left|\Pi_{\sigma \tau}^*\right| \leq \prod_{j=1,j\neq t}^m(\frac{s_{jt}^{2n_j}}{s_j^{n_j}}\prod_{p} p^{[\log(R+S(N+n_j))/\log p]}){n_t \choose \kappa}\kappa\prod_{p} p^{[\log(R+S(N+n_t))/\log p]} 
\]
\[
\leq (4S^3)^{N}e^{6(Rm+S(m+1)N)}.
\]
Next we combine this estimate, the upper bound
\[
\left|\prod_{\mu =0}^{N-1}\frac{\gamma_\sigma - \mu}{1+\mu}\right| \leq \prod_{p}p^{[\log(R+2NS)/\log p]} \leq e^{6(R+2SN)} 
\]
obtained by Lemma 2, and (\ref{8}) to obtain
\begin{equation}\label{16}
\left|k!a_k\right| \leq 2^k(8S^3)^Ne^{6(R(m+1)+S(m+3)N)} =: 2^kF_1^*.
\end{equation} 
An analog of (\ref{13}) is now
\begin{equation}\label{17}
\left|Q_0(z)\right| \leq \sum_{k=0}^N \left|a_{\sigma k}z^k\right| \leq F_1^*e^{2\left|z\right|}.
\end{equation}

The denominators of the coefficients $c_{0j\mu}$ of the polynomials $P_{0j}(z)$ can be considered similarly as the coefficients of $P_{ij}(z) \ (i = 1,\ldots,m)$ before, and these are factors of
\begin{equation}\label{18}
D_2^* := D_1^*\prod_{p} p^{[\log(R+SN)/\log p]} \leq S^{2(m+1)N}e^{6(m(\hat{R}+\hat{S}N)+2R+3SN))} =: E_2^*, 
\end{equation}
and clearly $D_2 \mid D_2^*$.

The above considerations lead to the following lemma.\\

\noindent
\textbf{Lemma 3}. \textit{Let $\alpha = a/b \neq 0$, where $a,b \in \mathbb{Z}_K$. Then
\[
\left|Q_i(\alpha)\right| \leq e^{c_1+c_2N}, \quad i = 0,1,\ldots,m,
\]
where
\[
c_1 = 6R(m+1)+2\left|\alpha\right|, \ c_2 = \log 8+3\log S+6S(m+3).
\]
Further, there exists an integer $D(N) \in \mathbb{Z}_K\setminus\{0\}$ such that
\[
(N+1)!D(N)Q_i(\alpha), \ (N+1)!D(N)P_{ij}(\alpha) \in \mathbb{Z}_K, \ i=0,1,\ldots,m;j=1,\ldots,m,
\]
and
\[
\left|D(N)\right| \leq e^{c_3+c_4N},
\]
where}
\[
c_3 = \log \left|b\right|+12R(1+Sm), \ c_4 = \log \left|b\right|+2(m+1)\log S+6S(3+Sm).\\
\]

\section{Remainder terms}

In this section we give an upper bound for the remainder terms.\\

\noindent
\textbf{Lemma 4}. \textit{We have
\[
\left|(N+1)!D(N)R_{ij}(\alpha)\right| \leq e^{c_5+c_6N}N^{-n_j}, \ i=0,1,\ldots,m; j = 1,\ldots,m,
\]
where}
\[
c_5 = c_1+c_3+\log 2+2(S^2-1)\left|\alpha\right|, \ c_6 = c_2+c_4+3\log 2+4\log S+2\log \max\{1,\left|\alpha\right|\}.\\
\]

\textit{Proof.} We first consider
\[
R_{0j}(z) = \sum_{\nu=N+n_j+1}^\infty c_{0j\nu}z^\nu,
\]
where, by (\ref{16}) and Lemma 2,
\[
\left|c_{0j\nu}\right| \leq \left|\sum_{k=0}^N\frac{k!a_k}{k![\nu-k]_j}\right| \leq F_1^*\sum_{k=0}^N\frac{2^k}{k!(\nu-k)!}\frac{(\nu-k)!}{\left|(\lambda_j+1)_{\nu-k}\right|} \leq 2^{N+1}F_1^*\frac{(2S^2)^\nu}{\nu!}.
\] 
Thus
\[
\left|R_{0j}(\alpha)\right| \leq 2^{N+1}F_1^*\sum_{\nu=N+n_j+1}^\infty \frac{\left|2S^2\alpha\right|^\nu}{\nu!} \leq 2^{N+1}F_1^*\frac{\left|2S^2\alpha\right|^{N+n_j+1}}{(N+n_j+1)!}e^{2S^2\left|\alpha\right|},
\]
and so, by (\ref{16}) and Lemma 3,
\[
\left|(N+1)!D(N)R_{0j}(\alpha)\right| \leq e^{c_5+c_6N}N^{-n_j}.
\]

For the consideration of $R_{ij}(z) \ (i\geq 1)$ we only need to replace above $F_1^*$ by $F_1$. This proves Lemma 4.\\

\section{Proof of Theorem 1}

Let us denote
\[
Q_i := (N+1)!D(N)Q_i(\alpha), \ P_{ij} := (N+1)!D(N)P_{ij}(\alpha), \ i =0,1,\ldots,m; j=1,\ldots,m.
\] 
By Lemma 3 all these numbers are integers in $K$, and
\[
\left|Q_i\right| \leq e^{N\log N+\hat{b}_1N+b_3}, \quad \hat{b}_1 = c_2+c_4+1, \ b_3 = c_1+c_3.
\]
Lemma 1 implies that the determinant
\[
\det (Q_i \ P_{i1} \ \ldots \ P_{im})_{i=0,1,\ldots,m} \neq 0.
\]
Further, by using Lemma 4, we see that if
\[
R_{ij} = Q_i \varphi_j(\alpha) - P_{ij}, \quad i=0,1,\ldots,m; j=1,\ldots,m,
\]
then
\[
\left|R_{ij}\right| \leq e^{-n_j\log N+\hat{e}_1N +e_3}, \quad \hat{e}_1 = c_6, \ e_3 = c_5.
\]
By denoting $b_1 = \hat{b}_1+1, e_1 = \hat{e}_1+1$, we have
\[
\left|Q_i\right| \leq e^{N\log N +b_1N}, \ \left|R_{ij}\right| \leq e^{-n_j\log N +e_1N}
\]
for all $i, j$ and $N \geq N_2 := \max \{b_3,e_3\}$.

The application of \cite[Corollary 3.5]{Ma} gives now the following result for linear forms $\Lambda = \beta_0+\beta_1\varphi_1(\alpha)+\cdots+\beta_m\varphi_m(\alpha)$, where $\alpha = a/b \in K\setminus \{0\}$ and $a,b,\beta_j \in \mathbb{Z}_K, (\beta_0,\beta_1,\ldots,\beta_m) \neq \underline{0}$. Let $\hat{H} = (2m)^mH, H = \prod_{j=1}^m h_j, h_j = \max \{1,\left|\beta_j\right|\} \ (j=1,\ldots,m)$, and let $x_2 = \max \{1,x\}$, where $x$ is the largest solution of the equation $x\log x = 2e_1m(x+m)$. If $(\beta_0,\beta_1,\ldots,\beta_m) \in \mathbb{Z}_K^{m+1}\setminus\{\underline{0}\}$ satisfies
\[
2\log \hat{H} \geq \max \{2\log N_2, x_2\log x_2, e^e\},
\]
then
\[
\left|\Lambda\right| > \frac{1}{2^{m+1}e^{m(1+b_1+e_1m)}}(\frac{\log\log \hat{H}}{\log \hat{H}})^m \hat{H}^{-1-\frac{4(1+b_1+e_1m)}{\log\log \hat{H}}}.
\]
Here 
\[
1+b_1+e_1m = d_0+d_1m+d_2m^2,
\]
where
\[
d_0=3+3\log 2+5\log S+36S+\log \left|b\right|,
\]
\[
\ d_1=1+6\log 2+9\log S+36S+\log \left|b\right|+2\log\max \{1,\left|\alpha\right|\},
\]
\[
d_2 = 2\log S+6S+6S^2,
\]
remember that $R = \max \{\left|r_j\right|\}$ and $S = \max \{s_j\}$, where $r_j/s_j = \lambda_j$. Thus
\[
\left|\Lambda\right| > H^{-1-\frac{6(d_0+d_1m+d_2m^2)}{\log\log H}}
\]
for all $H \geq H_0$, where $H_0$ is an effectively computable positive constant depending on $\lambda_1,\ldots, \lambda_m, m$ and $\alpha$. This proves Theorem 1.

\begin{small}

\vskip0.6cm
\noindent
Keijo V\"a\"an\"anen \\
Department of Mathematical Sciences  \\
University of Oulu \\
P. O. Box 3000 \\
90014 Oulun yliopisto, Finland \\
E-mail: keijo.vaananen@oulu.fi \\
\end{small}

\end{document}